# On the use of self-organizing maps to accelerate vector quantization

Eric de Bodt[1], Marie Cottrell[2], Patrick Letremy[2], Michel Verleysen[3]

[1] Université catholique de Louvain, IAG-FIN, 1 place des Doyens,
B-1348 Louvain-la-Neuve, Belgium
and
Université Lille 2, ESA, place Déliot, BP 381,
F-59020 Lille, France

[2] Université Paris I, SAMOS-MATISSE, 90 rue de Tolbiac,
F-75634 Paris Cedex 13, France

[3] Université catholique de Louvain, DICE, 3 place du Levant,
B-1348 Louvain-la-Neuve, Belgium

**Abstract.**
Self-organizing maps (SOM) are widely used for their topology preservation property: neighboring input vectors are quantified (or classified) either on the same location or on neighbor ones on a predefined grid. SOM are also widely used for their more classical vector quantization property. We show in this paper that using SOM instead of the more classical Simple Competitive Learning (SCL) algorithm drastically increases the speed of convergence of the vector quantization process. This fact is demonstrated through extensive simulations on artificial and real examples, with specific SOM (fixed and decreasing neighborhoods) and SCL algorithms.

## 1. Motivation

Vector quantization (VQ) is a widely used tool in many data analysis' fields. It consists in replacing a continuous distribution by a finite set of quantizers, while minimizing a predefined distortion criterion. Vector quantization may be used in clustering or classification tasks, where the aim is to determine groups (clusters) of data sharing common properties. It can also be used in data compression, where the aim is to replace the initial data by a finite set of quantified ones; labeling the quantified set and using the labels rather than the data themselves makes compression possible. Vector quantization is basically an unsupervised process. Supervised variants exist (LVQ1, LVQ2, in Kohonen [7]); in this last case, the *distortion criterion* takes class labels into account.

While most recent learning algorithms in various domains do not consider anymore the limitations due to computational load because of the increasing power of computers, this question is still important in VQ. The reason is that VQ is often used on huge databases in high-dimensional spaces; therefore the learning process of such VQ task may take hours or days to converge, especially if the number of clusters is large too.

In this paper, we study how, and up to which point, the Kohonen Self-Organizing Maps (SOMs) may be used to accelerate the VQ learning process. It is well known that SOMs accelerate vector quantization, compared to other, more traditional VQ algorithms. However, to our knowledge, it exists no quantitative experimental study on this topic. To fill this gap is the main purpose of this paper. After reminding formal definitions of the conventional VQ and SOM algorithms, and in particular of their distortion measures (section 2), we study the rate of convergence of conventional VQ and SOM algorithms, both on artificial databases where the exact solution of VQ can be computed, and on real ones (sections 3 and 5). Section 4 uses a hybrid method with a fixed-neighborhood SOM as initialization to a conventional VQ algorithm, in order to benefit both from the accelerated convergence of SOM and from the lower distortion error after convergence of conventional VQ; this hybrid algorithm is similar to the conventional SOM with decreasing neighborhood, and is aimed to better quantify the gain in convergence speed.

## 2. Vector quantization

### 2.1. Principle and distortion

Vector quantization consists in replacing a continuous distribution, in some cases known only through a finite number of samples, by a finite number of quantizers; the number of quantizers must be (much) smaller than the number of known samples. Each quantizer defines a cluster in the space; the principle of vector quantization is to project all samples in a cluster on the corresponding quantizer.

Most of the methods used to perform VQ necessitate setting a priori the number of clusters or quantizers. The choice of this number results from a trade-off between the precision (distortion) of the quantization and the necessity of an efficient description of the resulting clusters (quantity of information kept after quantization).

Once the number of quantizers is predefined, a good criterion of the quality of the clustering is the *distortion*, which measures the deviation between the data and the corresponding quantizers. Let us recall the main definitions and introduce our notations.

We consider a continuous data space $\Omega$, of dimension $d$, endowed by a continuous probability density function (pdf) $f(x)$, where the cumulated density (or repartition function) is $F(x)=P(X<x)$ (where $P$ is the probability law, and where the inequality is verified in each dimension).

A vector quantization $\Phi$ is an application from the continuous space $\Omega$ to a finite subset $F$ (the *codebook*) formed by $n$ *code-vectors* or *centroids* or *quantizers* $q_1, q_2, \ldots, q_n$ of $\Omega$. The positions of the code-vectors are supposed to be computed as a result of a *quantization algorithm* or *learning algorithm*.

The aim of a *vector quantization* (VQ) is to compress the information by replacing all elements $x$ of a cluster $C_i$ (subset of $\Omega$) by a unique quantizer (or code-vector, or centroid) $q_i$. For a given number $n$ of code-vectors, vector quantization tries to minimize the loss of information or *distortion*, measured by the mean quadratic error

$$\xi(f,\Phi) = \xi(f, q_1, q_2, \ldots, q_n) = \sum_{i=1}^{n} \int_{C_i} \left\| x - q_i \right\|^2 f(x)dx \,. \tag{1}$$

If a $N$-sample $x_1, x_2, \ldots, x_N$ is available (randomly chosen according to $f(x)$), this distortion is estimated by the intra-class sum of squares

$$\hat{\xi}(f,\Phi) = \hat{\xi}(f, q_1, q_2, \ldots, q_n) = \frac{1}{N} \sum_{i=1}^{n} \sum_{x_j \in C_i} \left\| x_j - q_i \right\|^2 \,. \tag{2}$$

All classical VQ algorithms (LBG, k-means,…) minimize this distortion function by choosing appropriate centroid locations. See for example Anderberg [1] or Bishop et al. [2] for proofs. There is no unique minimum of the distortion function, and the result strongly depends on the initialization.

## 2.2. Simple Competitive Learning and batch VQ algorithms

There exist many algorithms that deal with the VQ problem. Most of them are very slow in terms of convergence speed. The most popular one is the so-called Simple Competitive Learning algorithm (SCL) that can be defined as follows (see for example [5]):

let $\Omega$ be the data space (with dimension $d$), endowed with a density probability function $f(x)$. The data are randomly drawn according to the density $f(x)$ and are denoted by $x_1, x_2, \ldots, x_N$. The number of desired clusters is a priori fixed to be $n$. The quantizers $q_1, q_2, \ldots, q_n$ are randomly initialized. At each step $t$,

- a data $x_{t+1}$ is randomly drawn according to the density $f(x)$ ;
- the winning quantizer $q_{win(t)}$ is determined by minimizing the classical Euclidean norm
  $$\| x_{t+1} - q_{win(t)} \| = \min_i \| x_{t+1} - q_i \| \,;$$
- the quantizer $q_{win(t)}$ is updated by $q_{win(t+1)} = q_{win(t)} + \varepsilon(t)\,(x_{t+1} - q_{win(t)})$.

where $\varepsilon(t)$ is an adaptation parameter which satisfies the classical Robbins-Monro conditions ($\Sigma\, \varepsilon(t) = \infty$ and $\Sigma\, \varepsilon^2(t) < \infty$).

The SCL algorithm is in fact the *stochastic or on-line* version of the Forgy algorithm (also called *moving centers algorithm*, LBG or Lloyd's algorithm, see for example [4], [8], [9]). In that version of the algorithm, the quantizers are randomly initialized. At each step $t$, the clusters $C_1, C_2, \ldots, C_n$ are determined by putting in class $C_j$ the data which are closer to $q_j$ than to any other quantizer $q_i$. Then the mean values in each cluster are simultaneously computed and taken as new quantizers. The process is then repeated. The Forgy algorithm works off-line as a batch algorithm. It will be referred to as BVQ (for Batch VQ) in the following. It also exits an intermediate version of the algorithm, frequently named the K-means method [10]. In that case, at each step, only one data is randomly chosen, and only the winning quantizer is updated as the mean value of its class.

It can be proven ([1], [2]) and it is well-known that BVQ minimizes the *distortion* (1) or, more exactly, the estimated one (2). Note that the stochastic SCL algorithm also minimizes this distortion, but only in mean value: at each step, there is a positive probability to increase the distortion, as for any stochastic algorithm.

Let us denote by $q_1^*, q_2^*, \ldots, q_n^*$ one set of quantizers which (locally) minimizes the distortion. This minimum needs not to be unique and generally depends on the initial values. At a local minimum of the distortion, each $q_i^*$ is the center of gravity of its class $C_i$, with respect to the density $f$, and the quantizers are fixed points of the BVQ algorithm. In an exact form,

$$q^*_i = \frac{\int_{C_i} x\, f(x)\, dx}{\int_{C_i} f(x)\, dx}, \tag{3}$$

estimated by

$$\hat{q}^*_i = \frac{\sum_{x_j \in C_i} x_j}{\sum_{x_j \in C_i} 1}. \tag{4}$$

Note that the equations are implicit ones, since the $C_i$ are defined according to the positions of the $q^*_i$.

For example, in the one-dimensional case, the classes $C_i$ ($1 \leq i \leq n$) are intervals defined by $C_i = [a_i, b_i]$, with $a_i = ½\,(q^*_{i-1} + q^*_i)$ and $b_i = ½\,(q^*_{i+1} + q^*_i)$, for $1 < i < n$, and $a_1 = \inf(\Omega)$, $b_n = \sup(\Omega)$.

The main goal of this paper (after a preliminary work presented to ESANN'99, see [3]) is to evaluate the speed of convergence of VQ algorithms. In situations where the solution is unique and where it is possible to compute the exact values $q^*_i$, the performances will be evaluated through the rate at which the values $q_i$ converge to $q^*_i$ (see section 3).

## 2.3. Self-Organizing Maps (SOM)

Let us consider now the SOM algorithm (as defined by T. Kohonen in [6]). As we will show it, it can be seen as an extension of the Simple Competitive Learning Algorithm (SCL) in its classical stochastic form, and of the Forgy Algorithm (BVQ) in its batch form. We will consider here the SOM algorithm with a fixed number of neighbors (although the number of neighbors uses to decrease with time in practical implementations).

For a given neighborhood structure, where $V(i)$ denotes the neighborhood of unit $i$, the SOM algorithm is defined as follows. The quantizers $q_1, q_2, \ldots, q_n$ are randomly initialized. At each step $t$,

- a data $x_{t+1}$ is randomly drawn according to the density $f(x)$ ;
- the winning quantizer $q_{win(t)}$ is determined by minimizing the classical Euclidean norm

$$\| x_{t+1} - q_{win(t)} \| = \min_i \| x_{t+1} - q_i \| ;$$

- the quantizer $q_{win(t)}$ and its neighbors $q_k(t)$ for $k$ in $V(j)$ are updated by

$$q_{win,k(t+1)} = q_{win,k(t)} + \varepsilon(t)\,(x_{t+1} - q_{win,k(t)}).$$

where $\varepsilon(t)$ is an adaptation parameter which satisfies the classical Robbins-Monro conditions ($\sum \varepsilon(t) = \infty$ and $\sum \varepsilon^2(t) < \infty$).

We see that the SCL algorithm is a particular case of the SOM algorithm, when the neighborhood is reduced to zero. Sometimes SCL is called *0-neighbor Kohonen algorithm*. There also exists a batch SOM algorithm, similar to the Forgy (BVQ) algorithm. The only difference is that at each step, for a given set of classes $C_1, C_2, \ldots, C_n$, the quantizer $q_j$ is set to the mean value of the union of the class $C_j$ and its neighbors. See [7] for example.

It appears clearly that the SOM algorithm is different from the SCL algorithm only because a neighborhood structure is defined between the $n$ quantizers. The neighborhood structure of the SOM algorithm is most frequently used for visualization and data interpretation properties. We however only consider it here as a way to accelerate the convergence of the SOM algorithm. In other words, we are only interested here in the VQ property of self-organizing maps, and not in its topological properties.

In the one-dimensional case, and for a one-dimensional structure of neighborhood, if the neighborhood $V(i)$ contains $i$-1, $i$, $i$+1 (two-neighbor case), the limit points $q^*_i$ of the batch SOM algorithm verify equation (3) or equation (4), where $C_i$ is replaced by $C_i^2 = C_{i-1} \cup C_i \cup C_{i+1} = [a_i, b_i]$, with $a_i = ½ (q^*_{i-2} + q^*_{i-1})$ and $b_i = ½ (q^*_{i+1} + q^*_{i+2})$, for $2 < i < n-1$, and $a_1 = a_2 = \inf(\Omega)$, and $b_{n-1} = b_n = \sup(\Omega)$.

Here again, the batch SOM algorithm is nothing else than the iterative computation of the solutions of equations (3) or (4), when $C_i$ is replaced by $C_i^2$.

The batch SOM algorithm and the classical stochastic SOM algorithm do not decrease anymore the distortion (1), but a generalized distortion that is the distortion extended to the neighbor classes (as long as the number of neighbors $\nu$ remains fixed) [11]. This generalized distortion is given by

$$\xi_\nu(f, \Phi) = \xi_\nu(f, q_1, q_2, \ldots, q_n) = \sum_{i=1}^{n} \int_{\bigcup_{k \in V(i)} C_k} \|x - q_i\|^2 f(x)dx\,, \tag{5}$$

where $V(i)$ is the set of indexes in the neighborhood of $i$, including $i$. This generalized distortion function can also be estimated through a finite set of samples $x_1, x_2, \ldots, x_N$, similarly to (2).

## 3. Experimental results: convergence to the exact solution of the VQ problem

The SOM algorithm is not equivalent to the SCL algorithm: it is deemed to minimize the generalized distortion (5), and not the classical distortion of VQ problems (1). Despite this fact, we will show that the SOM algorithms perform better than the classical SCL algorithm, i.e. converge faster towards the solution of (1), at least during the first iterations of the algorithm.

We study this phenomenon from two points of view. First, in some cases where it is possible to exactly compute the solutions of equations (3) or (4), we evaluate the error between the current values and the optimal values as a function of the number of iterations, for both the SCL and SOM algorithms. This is the topic of this section. Secondly, for more realistic data, we compare the decreasing slope of the true distortion (1) as a function of the number of iterations, also for both algorithms. This is the topic of section 5

In some one-dimensional cases ($d$ = 1), if the set $\Omega$ is a real interval, and if the density $f$ is known and well-behaved, it is possible to directly compute the solutions $q^*_i$ , starting from a given set of increasing initial values, by an iterative numerical procedure.

If the initial values are ordered, the current values $q_1, q_2, \ldots, q_n$ will remain ordered at each iteration of the SCL algorithm. As mentioned in the previous section, the classes $C_i$ ($1 \le i \le n$) are therefore intervals defined by $C_i = [a_i, b_i]$, with $a_i = ½ (q_{i-1} + q_i)$ and $b_i = ½ (q_{i+1} + q_i)$, for $1 < i < n$, and $a_1 = \inf(\Omega)$, $b_n = \sup(\Omega)$. This constitutes the first set of equations ($C_i$ as a function of $q_i$) used in this iterative procedure.

Equations (3) or (4) have no explicit solutions. However, it is possible to compute analytically the solutions $q_i$, as a function of the limits $a_i$ and $b_i$ of the intervals $C_i$, for some "easy" densities $f(x)$. This will constitute the second set of equations (equivalent to (3)) used in the iterative procedure. Table 1 presents these recurrent explicit equations for the densities $f(x) = 2x, 3x^2, e^{-x}$.

| Density $f$ | Distribution Function | $a_0$ | $b_0$ | $q_i$ |
|---|---|---|---|---|
| $2x$ on [0,1] | $x^2$ | 0 | 1 | $q_i = \frac{2}{3}\frac{b_i^3 - a_i^3}{b_i^2 - a_i^2}$ |
| $3x^2$ on [0,1] | $x^3$ | 0 | 1 | $q_i = \frac{3}{4}\frac{b_i^4 - a_i^4}{b_i^3 - a_i^3}$ |
| $e^{-x}$ on [0,+∞] | $1 - e^{-x}$ | 0 | +∞ | $q_i = \frac{a_i e^{-a_i} + e^{-a_i} - b_i e^{-b_i} - e^{-b_i}}{e^{-a_i} - e^{-b_i}}$ |

**Table 1:** Exact computation of the quantizers as a function of the limits of the clusters, for some "easy" examples of densities.

Using alternatively the two sets of equations leads to the convergence to the optimal values $q^*_i$ of the quantizers. This iterative procedure is similar to the BVQ algorithm. Formulas in Table 1 are the analytical solutions of equation (3), while the BVQ algorithm usually involves in practical experiments the use of approximation (4).

Knowing the optimal values of the quantizers, it is possible to study the speed of convergence of any Vector Quantization algorithm. In that purpose, we will study the Euclidean distance between the current values of $(q_i(t))$ resulting from some VQ algorithm and the solutions $(q^*_i)$, as a function of the numbers of iterations. We define the mean quadratic error

$$D^2(t) = D^2(q(t),q^*) = (1/n)\, \Sigma_{1\leq i\leq n}\, (q_i(t) - q^*_i)^2 \qquad (6)$$

which will be the *error measure* of the Vector Quantization algorithm into consideration.

In practical situations, one can observe that the error measure $D^2(t)$ decreases to 0 very slowly when using the SCL algorithm. Note that in all simulations with several algorithms we carefully started from the same initial increasing values $q_i(0)$, including for the exact computation of the $(q^*_i)$, in order to avoid any effect due to the initial conditions.

In Figures 1, 2 and 3, we represent the variations of the error measure $D^2(t)$, for the SCL algorithm and for the SOM with 2, 4 and 8 neighbors; 100 quantizers were used for all simulations on artificial data. Figures 1, 2 and 3 respectively concern densities $f(x) = 2x$, $3x^2$ and $e^{-x}$. We can see for example that the SOM with neighbors decreases to the optimal values $(q^*_i)$ much faster than the SCL algorithm, even if it finally converges to its own optimal points. These optimal points minimize the generalized distortion extended to the neighbors (5), and are different from the $(q^*_i)$ (see section 2).

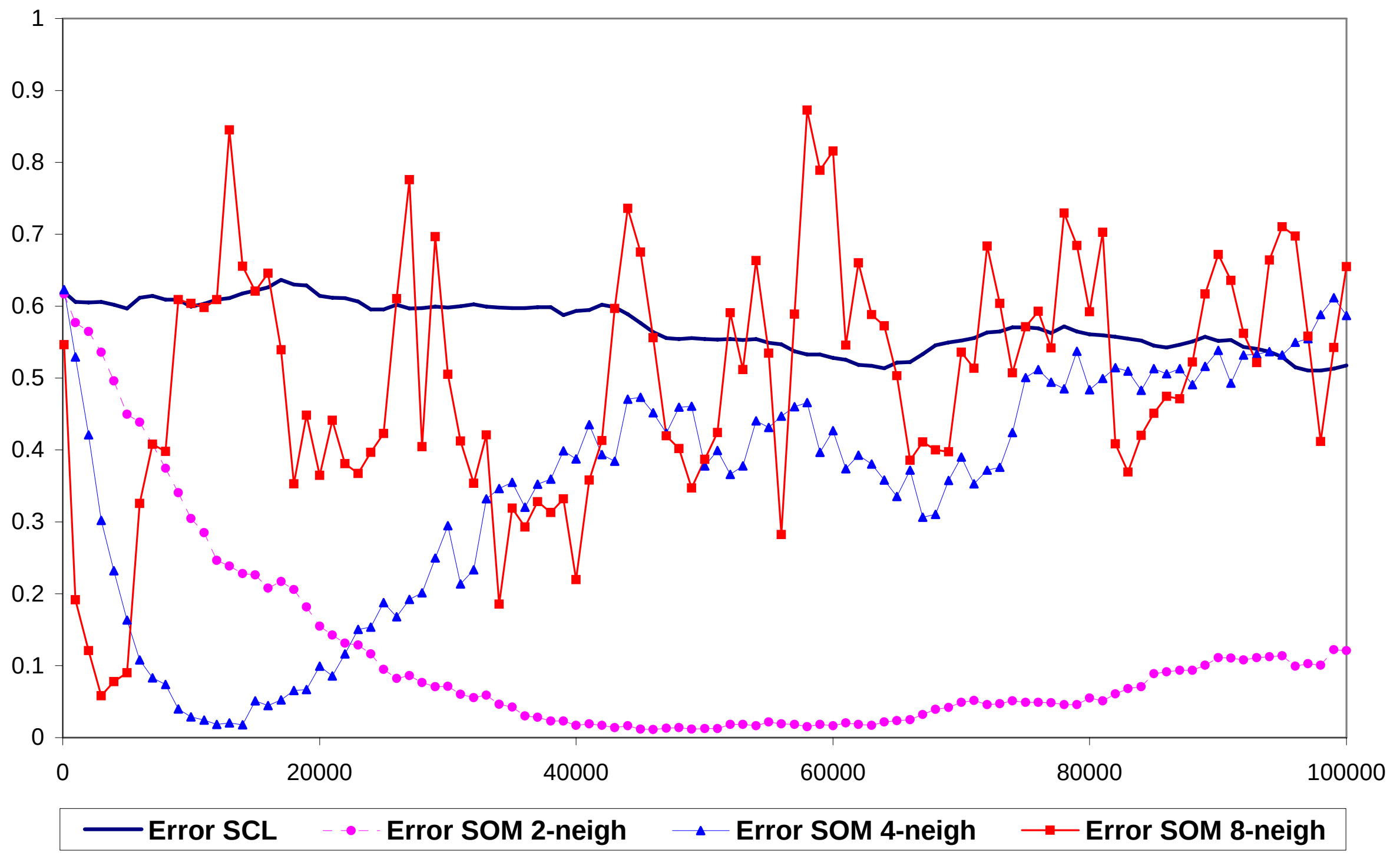

**Figure 1:** Evolution of $D^2(t)$ as a function of the number of iterations, for the density $f(x) = 2x$.

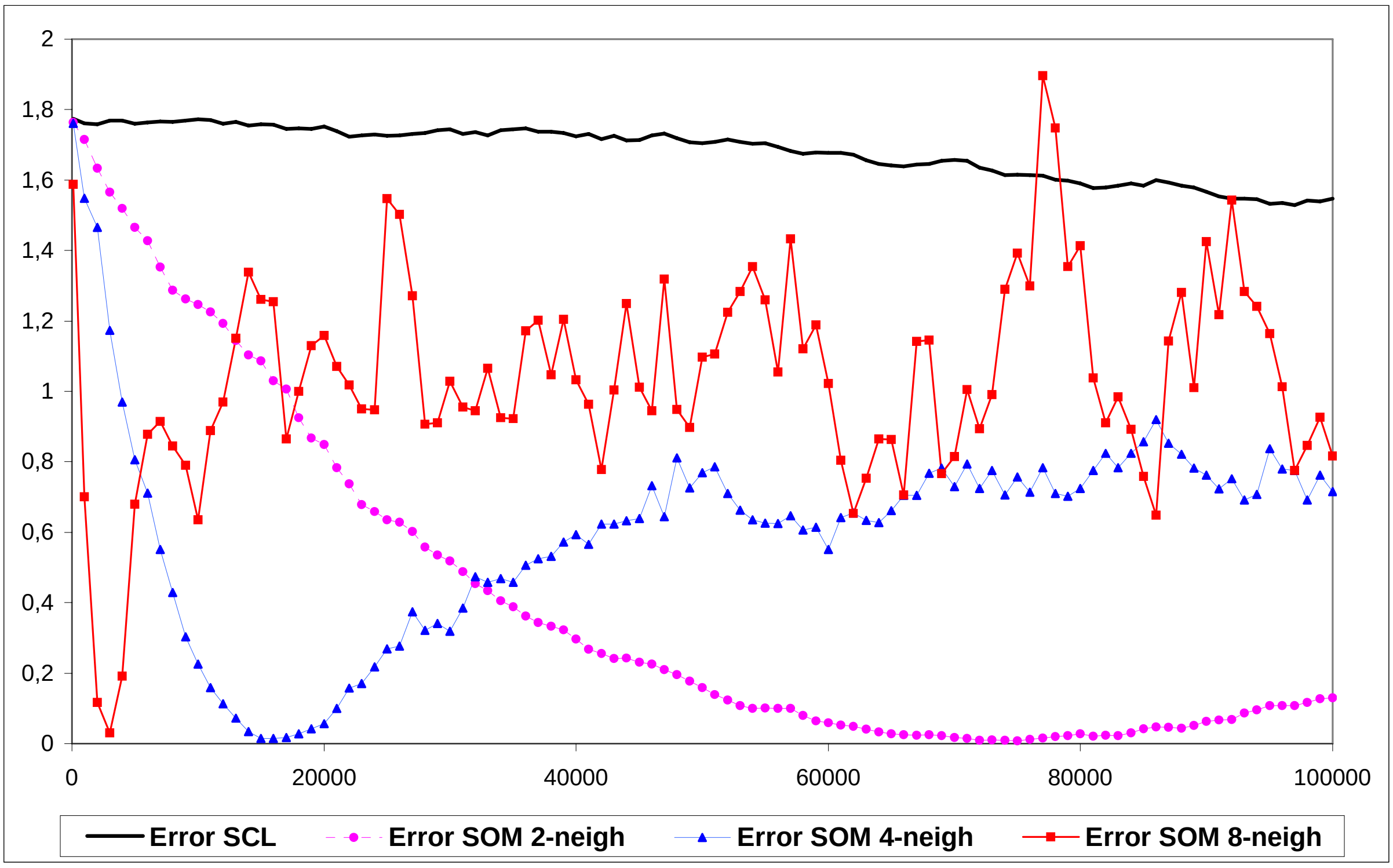

**Figure 2:** Evolution of $D^2(t)$ as a function of the number of iterations, for the density $f(x) = 3x^2$

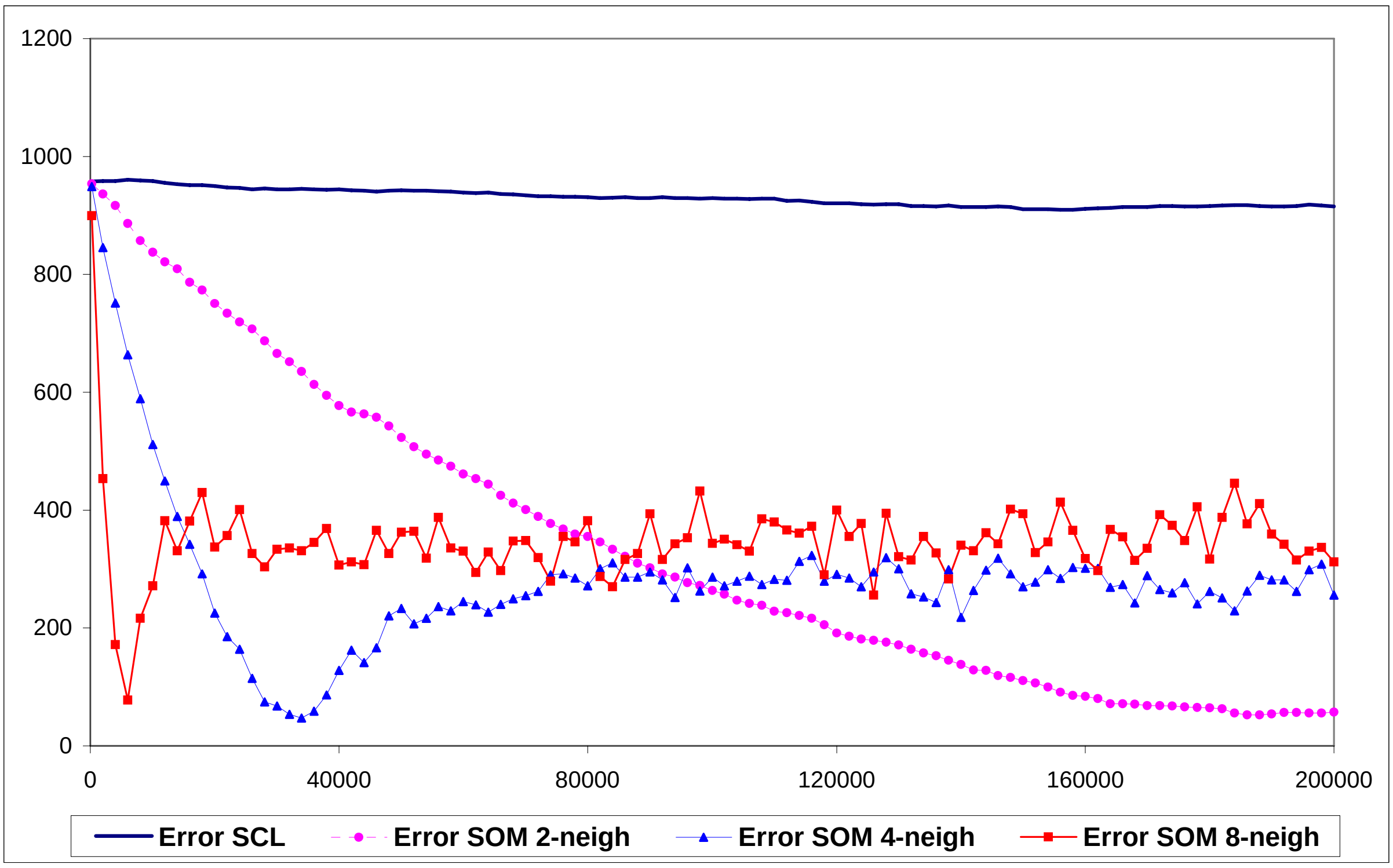

**Figure 3:** Evolution of $D^2(t)$ as a function of the number of iterations, for the density $f(x) = e^{-x}$

We also measured the evolution of $D^2(t)$ as a function of the number of iterations, for a Gaussian density $N(0,1)$. In this case, the exact values $q^*_i$ have been obtained through equation (4), by using very large samples to compute at each step the $C_i$. Figure 4 shows this evolution of $D^2(t)$ as a function of the number of iterations, respectively for the SCL algorithm and for the SOM with 2, 4, 8 and 16 neighbors.

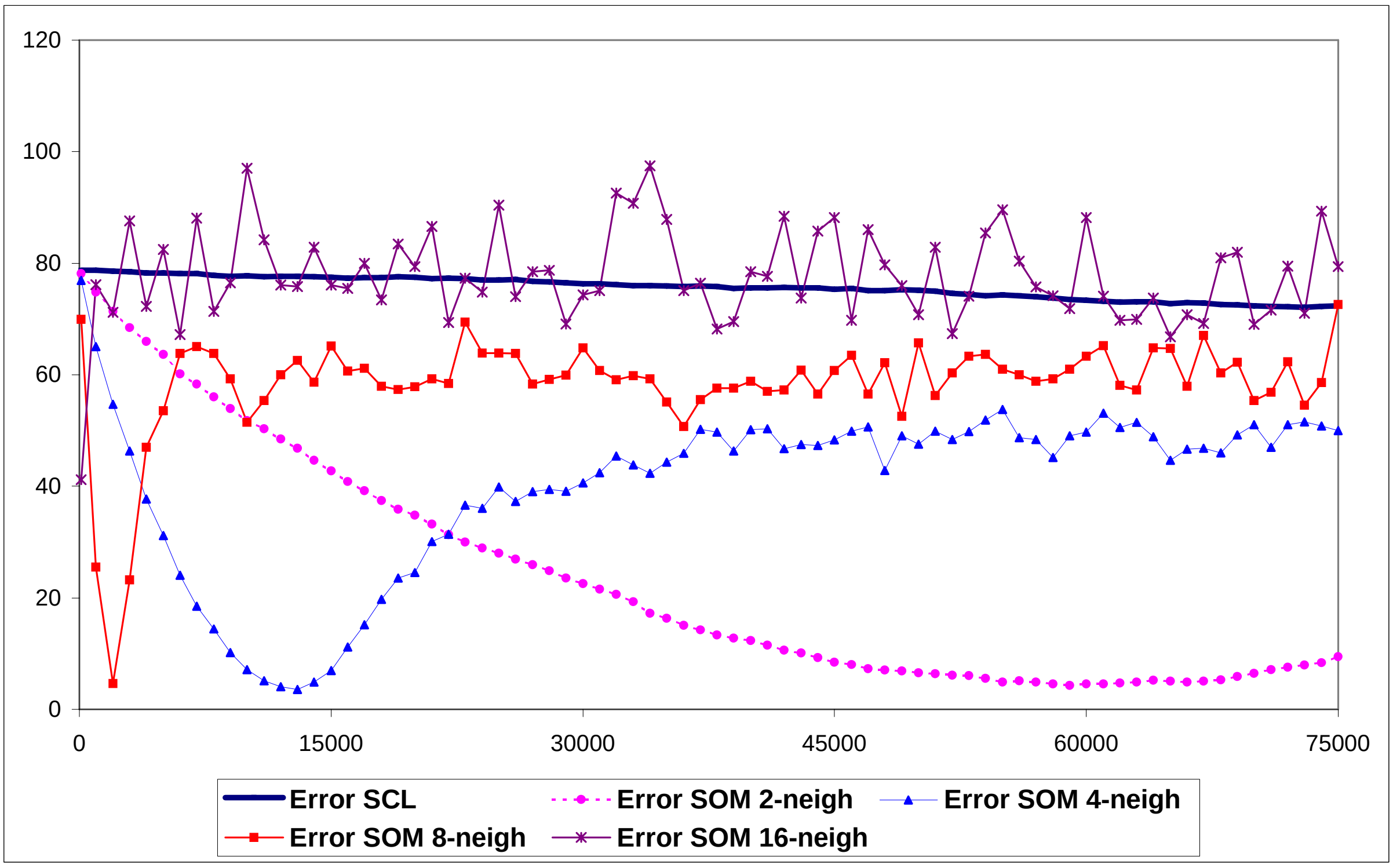

**Figure 4:** Evolution of $D^2(t)$ as a function of the number of iterations, for the standard Gaussian $N(0,1)$ density

One could argue that the comparisons are made on different algorithms, where the processing time per iteration is different. The comparison in terms of number of iterations is thus not fair if the total processing time is searched for. Nevertheless, as an example, the difference of the processing time in our simulations of a single iteration when using the 2-neighbors SOM algorithm instead of the SCL algorithm is significantly less than 1%. The differences shown in Figures 1 to 4 thus remain striking.

This section has shown that the use of SOM can greatly increase the speed of convergence towards the exact solutions of the VQ problem. Nevertheless, it must not be forgotten that the SOM algorithm will not finally converge to these solutions but rather to a minimum of (5). In the next section, we therefore use a mixed algorithm, beginning by some iterations of the SOM algorithm and ending with a classical SCL procedure, in order to benefit both from the accelerated convergence and from the convergence towards optimal states.

## 4. Hybrid algorithm SOM/SCL

Based on the results of the previous section, we propose to use a hybrid VQ algorithm (denoted by KSCL for Kohonen SCL), which consists of an initial phase (a SOM algorithm with $\nu$ neighbors), followed by the classical SCL. We compare the value of the error $D^2(t)$ after the same number of iterations for KSCL and SCL. Note that KSCL is not far from the classical SOM algorithm with decreasing neighborhood; the reason to use KSCL is not to suggest another algorithm, but to better quantify the gain in convergence speed compared to classical SCL, regardless of the specific decreasing function used in a classical SOM.

For example, let us fix a total number of iterations $T$, the initial ordered points $q_1(0)$, $q_2(0)$, …, $q_n(0)$, a constant $\varepsilon$ and several probability functions : $f(x) = 2x$ on [0,1], $f(x) = 3x^2$ on [0,1], $f(x) = e^{-x}$ on $[0, +\infty[$, and the standard Gaussian N(0,1). Let us also consider the 2-neighbors SOM algorithm ($\nu = 2$).

In Figures 5, 6, 7 and 8, we represent the evolution of the error measure for different KSCL algorithms and for the four probability densities that we took as examples. We consider four KSCL variants where the 2-neighbors SOM algorithm is used respectively during 0%, 30%, 60%, 90% of the total number of iterations $T$.

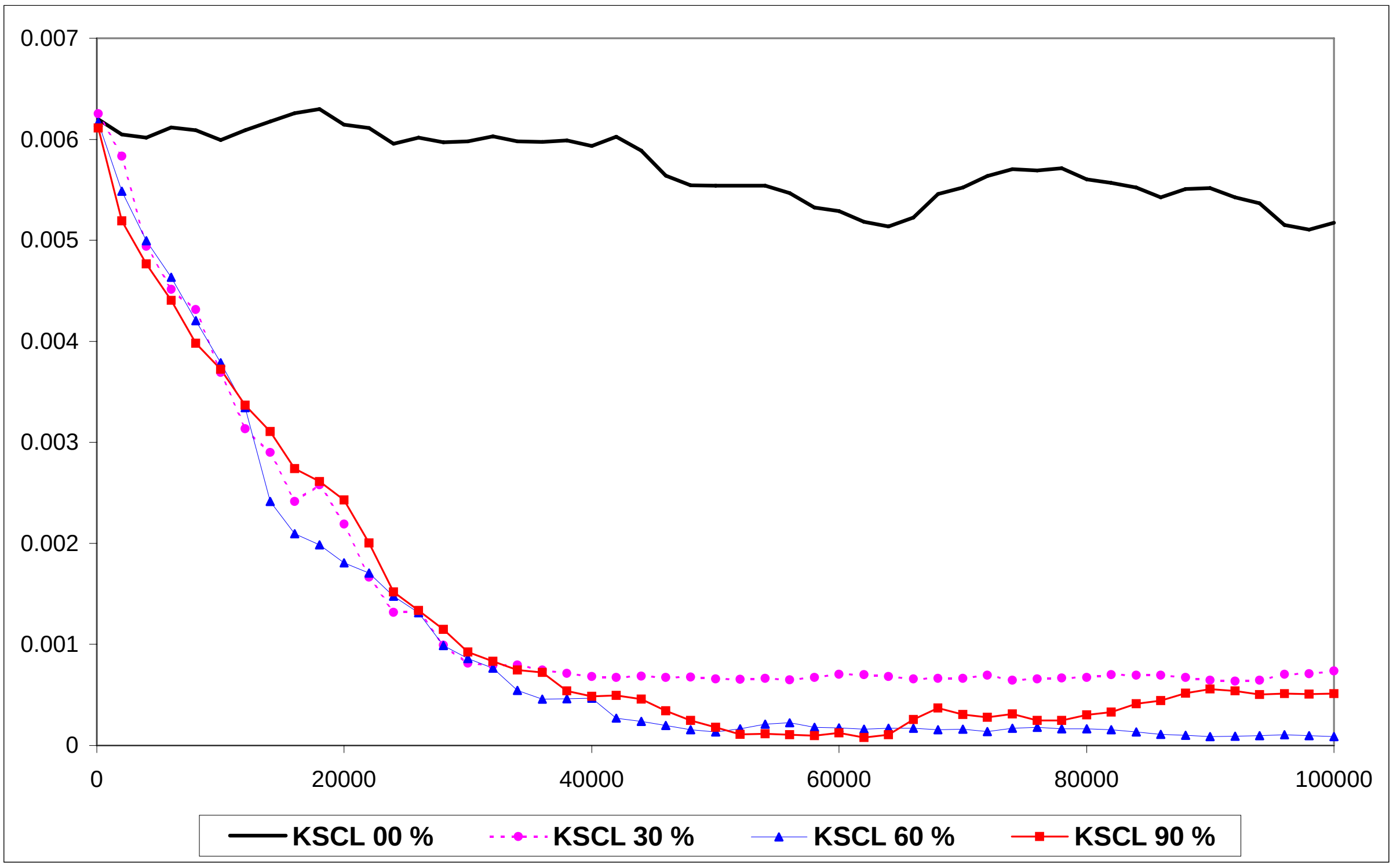

**Figure 5:** Evolution of $D^2(t)$ as a function of the number of iterations, for 4 variants of the KSCL algorithm, on the $2x$ density.

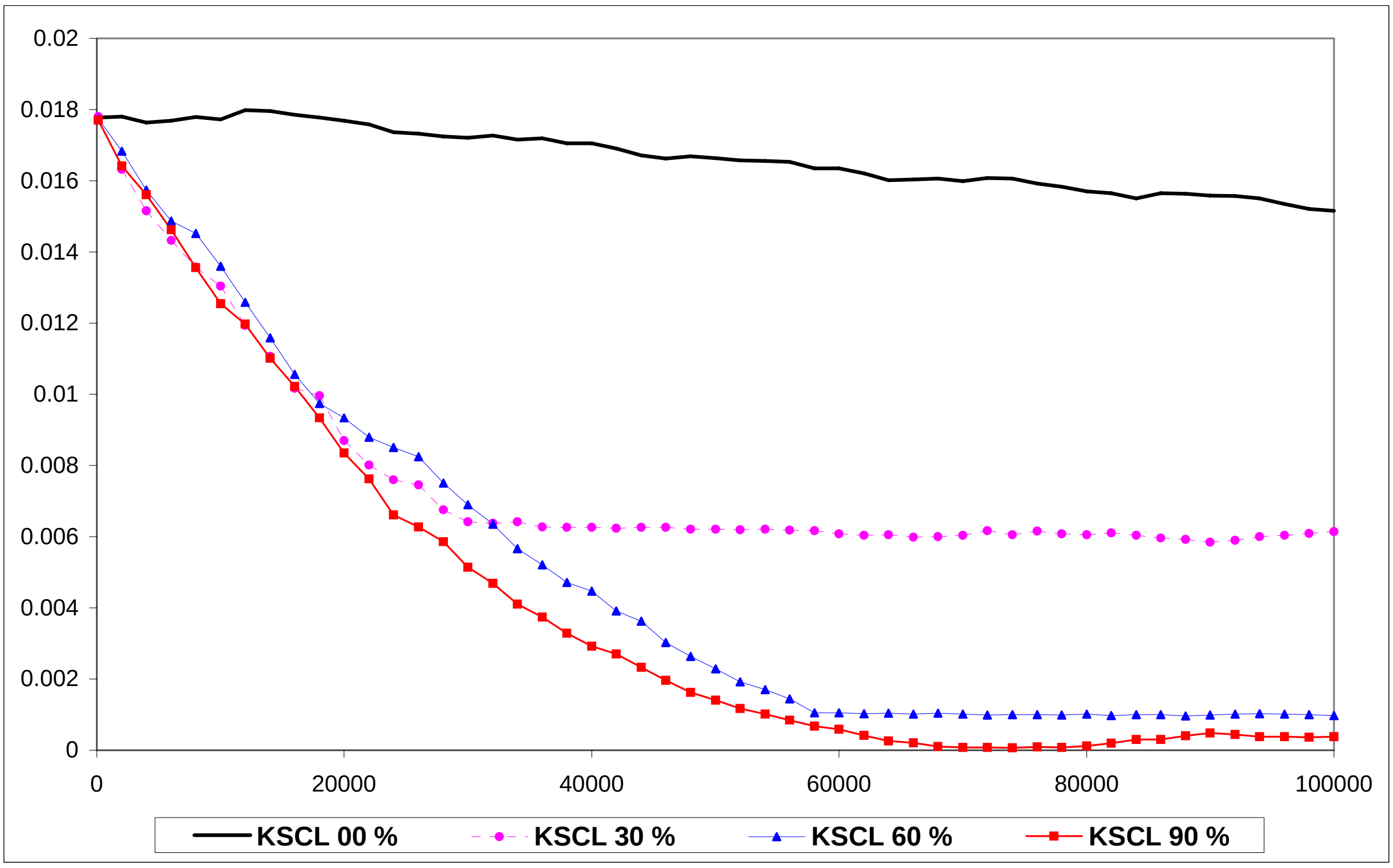

**Figure 6**: Evolution of $D^2(t)$ as a function of the number of iterations, for 4 variants of the KSCL algorithm, on the $3x^2$ density.

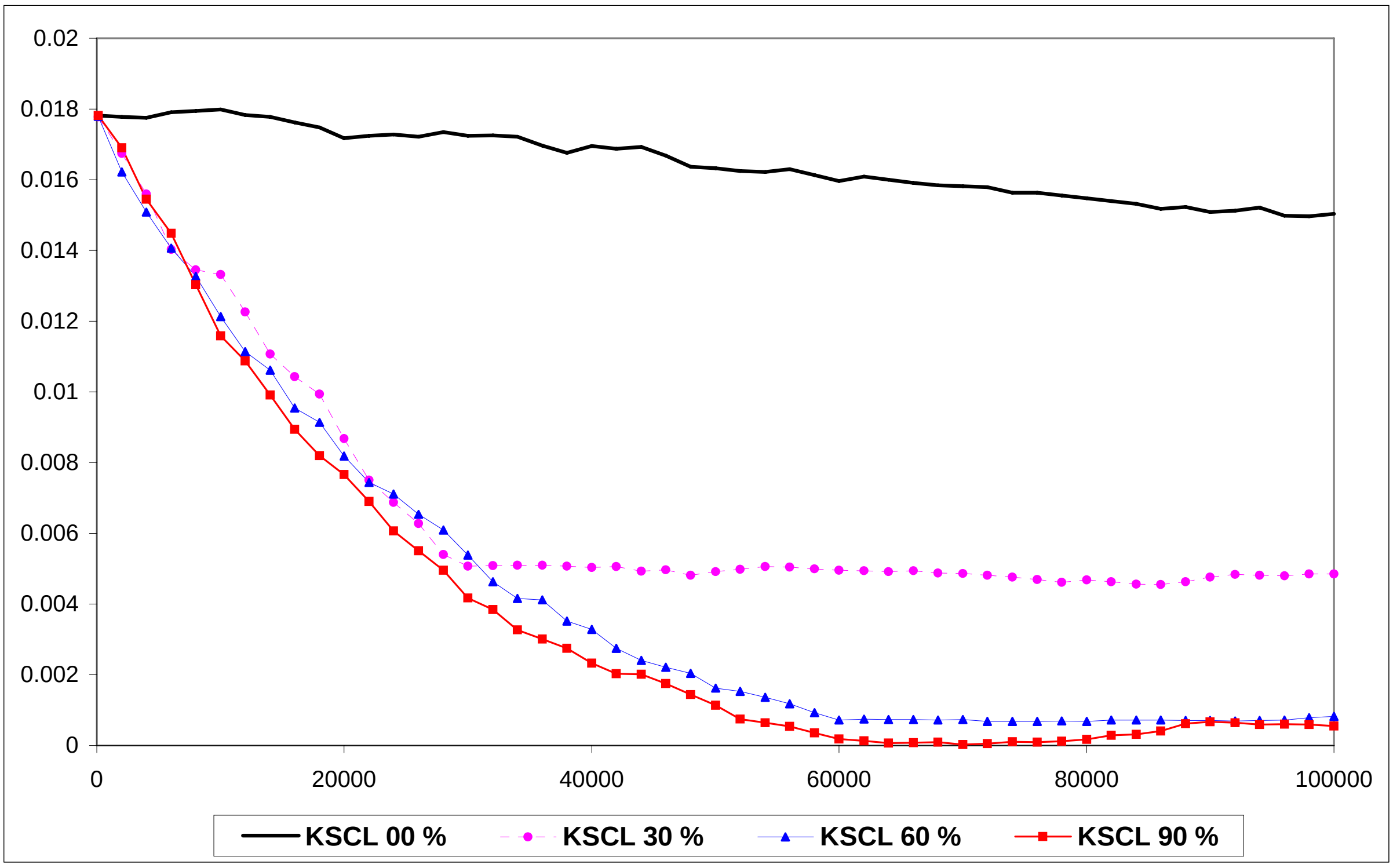

**Figure 7**: Evolution of $D^2(t)$ as a function of the number of iterations, for 4 variants of the KSCL algorithm, on the $e^{-x}$ density.

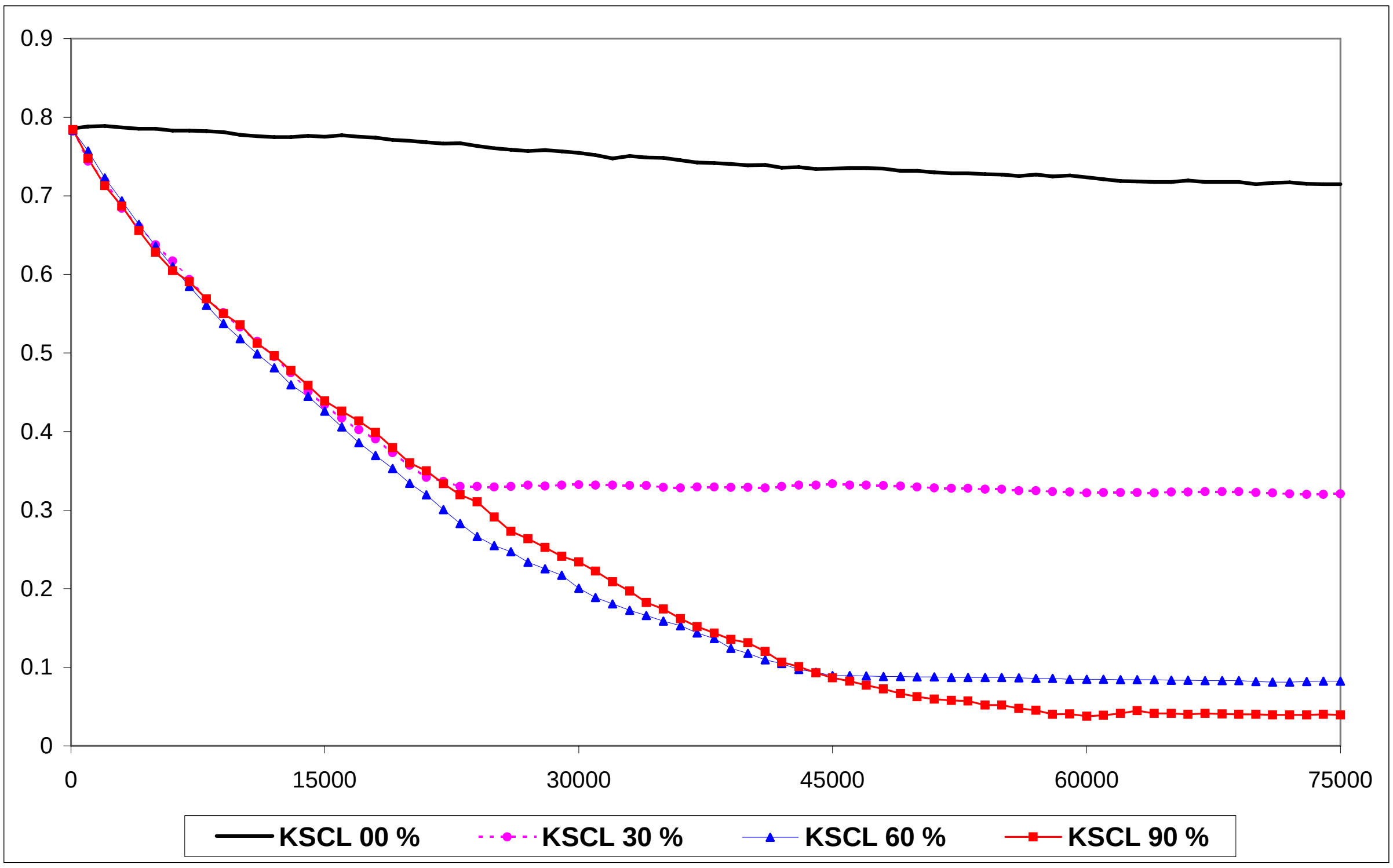

**Figure 8**: Evolution of $D^2(t)$ as a function of the number of iterations, for 4 variants of the KSCL algorithm, on the standard Gaussian density.

We can observe in all simulations that the 2-neighbors algorithm greatly accelerates the decrease of the error measure in the beginning of the curves. In all cases, using too early the SCL algorithm slows down the decrease.

Moreover, the performances remain better than those of the SCL algorithm, whatever is the choice of the KSCL variant. Nevertheless, it is also clear that determining the optimal iteration for substituting SOM by SCL strongly depends on the probability density. An optimal choice of this parameter would thus require extensive simulations, which is not the goal searched for here.

We may conclude this section by claiming that, in any case, *the SOM algorithm with a fixed number of neighbors can work as an efficient initialization of the SCL algorithm to accelerate the convergence and improve the performances.* We verified this statement on many other probability densities, real data and for several values of the number of neighbors' $v$.

In practical implementations of the SOM algorithm, the number of neighbors is made decreasing. Our observations confirm that this widely used strategy is very efficient to improve the decreasing of the *error measure* (6).

It would be interesting to consider this so-called error measure in multidimensional settings. Nevertheless, this concept is not well suited to dimensions greater than 1. The lack of ordering concept in dimension greater than 1 does not facilitate the problem and the correspondence between the current quantizers at a given iteration and their optimal values looses its clear meaning.

We will thus replace the concept of error measure by the concept of distortion as defined by (1) or (2). In the next section, we study how the distortion is decreasing along the quantization process, in both cases (SCL without neighbors, or SOM with a decreasing number of neighbors).

## 5. Experimental results: comparative evolution of the distortion

In this part, we study the vector quantization speed of SCL and of SOM, by computing the distortion defined in equation (2) in the case of real data.
As we mentioned previously, the SCL is supposed to minimize this distortion, while the SOM (with fixed or decreasing number of neighbors) is not. However we can observe that in any case, the SOM algorithm accelerates the decrease of the distortion, at least during a large part of the simulation.

We represent the distortion as a function of the number of iterations, for 5 different quantization algorithms: SCL and 4 variants of the SOM algorithm which differ by the number of neighbors. For a two-dimensional neighborhood structure, we consider successively 3 SOM algorithms with a fixed number of neighbors (SOM5, SOM9 and SOM25, the suffix being the number of neighbors) and then the classical SOM algorithm with a decreasing number of neighbors (from 25 to 1, the last part of the SOM iterations being equivalent to SCL).

We illustrate these simulations on two data sets. The table SAVING contains 5 variables measuring economic ratios for 42 countries between 1960 and 1970; the table TOP500 contains 6 variables relative to 77 American companies in 1986. [1]

Figure 9 represents the distortion for the SAVING data, with 25 quantizers, square grid (5 by 5) for the SOM algorithms and 500 steps of iterations. Figure 10 represents the distortion for the data SAVING, with 100 quantizers, square grid (10 by 10) for SOM algorithms and 1000 iterations. Figure 11 represents the distortion for the data TOP500, with 100 quantizers, square grid (10 by 10) for SOM algorithms and 1000 iterations.

[1] The data are available on the http://www.dice.ucl.ac.be/neural-nets/Research/Projects/projects.htm WEB page.

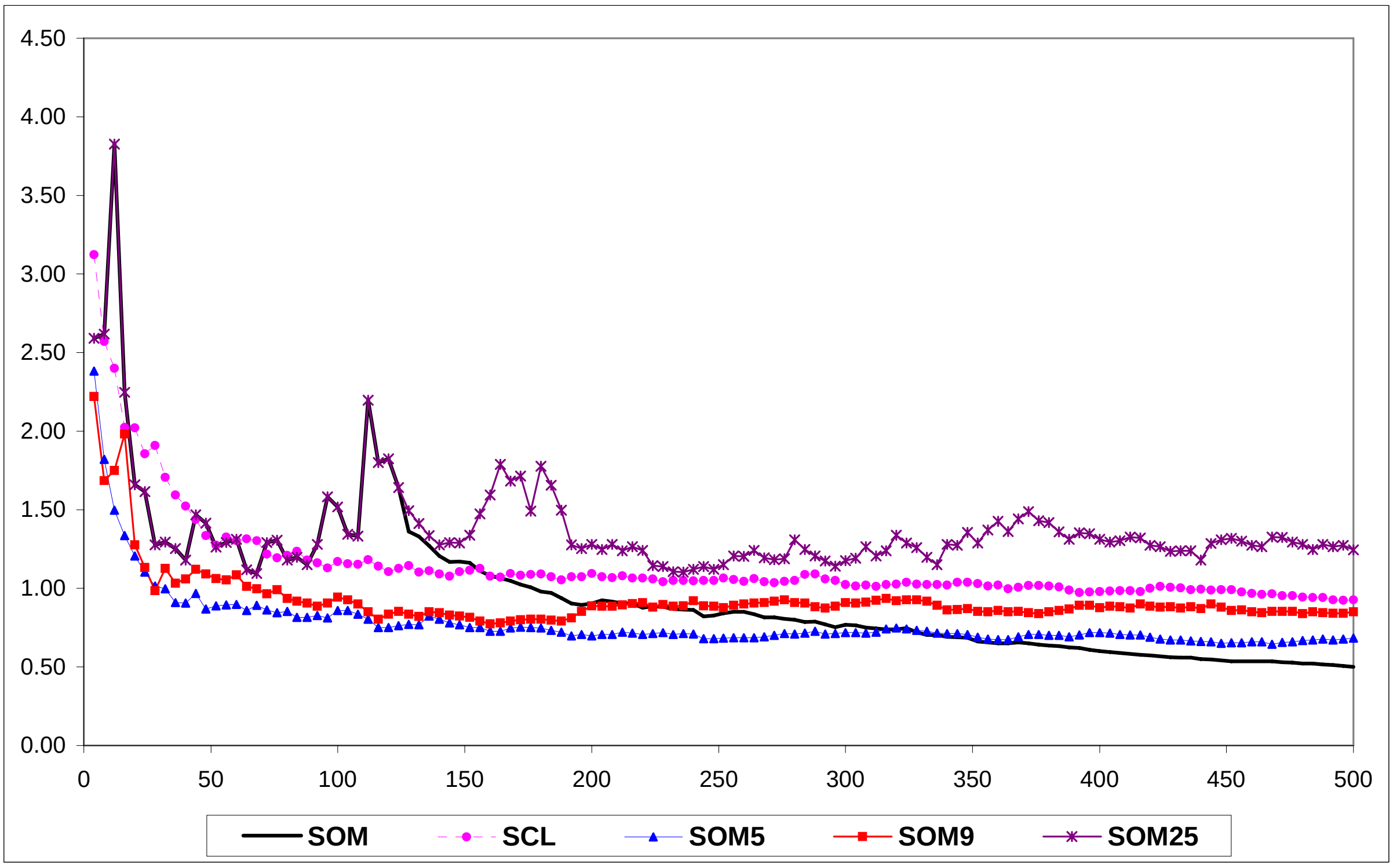

**Figure 9**: Evolution of the distortion as a function of the number of iterations, on the Saving data set; 25 quantizers are used (see text for details on the algorithms). The distortion is illustrated *after* each iteration.

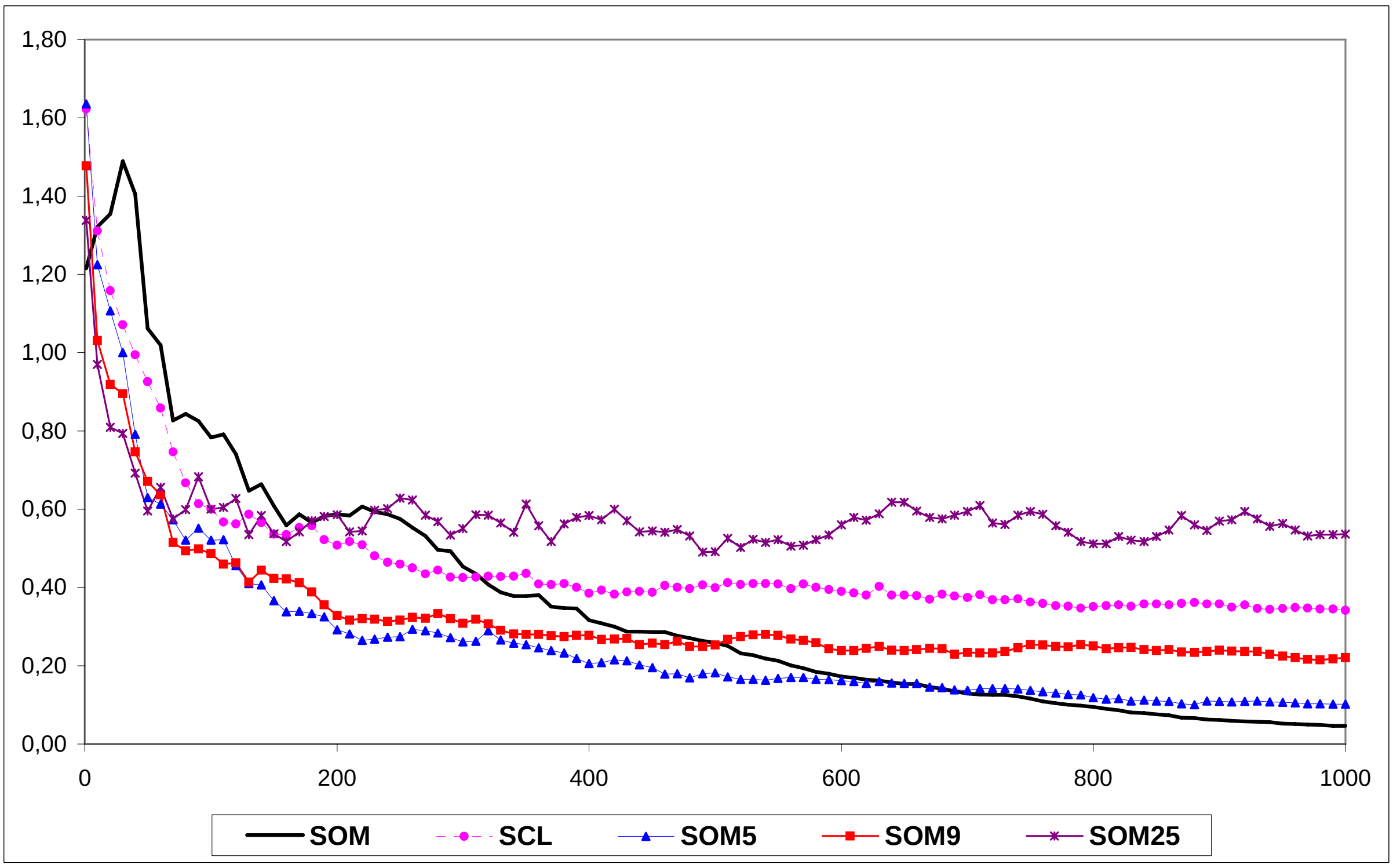

**Figure 10**: Evolution of the distortion as a function of the number of iterations, on the Saving data set; 100 quantizers are used (see text for details on the algorithms).

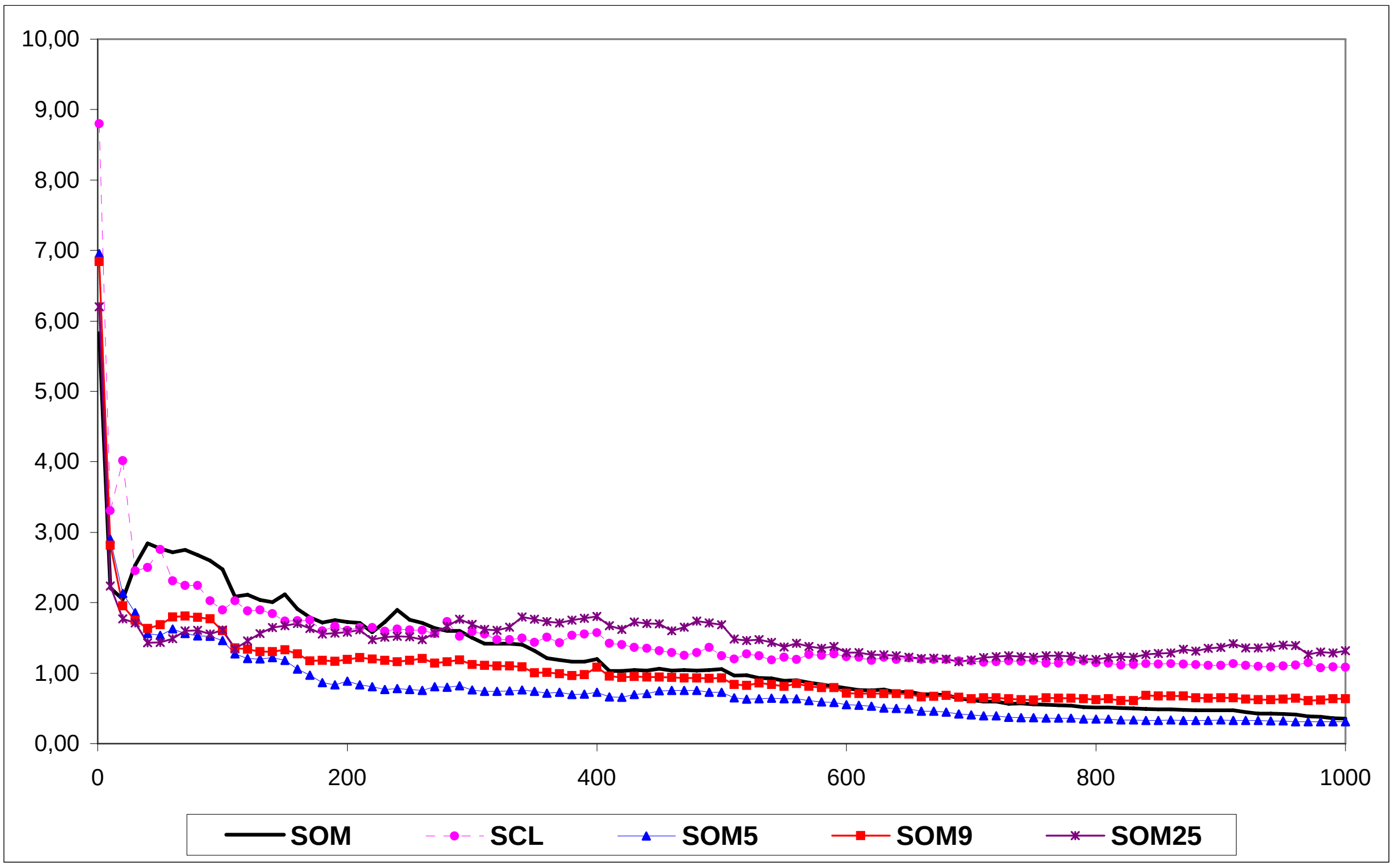

**Figure 11**: Evolution of the distortion as a function of the number of iterations, on the Top500 data set; 100 quantizers are used (see text for details on the algorithms).

In each simulation, we can see that the SOM algorithm performs as the best quantizer (it leads to a lower minimum of the distortion function). The SCL algorithm is very slow, and the SOM with non-decreasing number of neighbors is powerful at the beginning of the iterations, but allows at some iteration the distortion to increase (or would allow, after a larger number of iterations). In fact the classical SOM algorithm ending with no neighbor appears to be an excellent VQ algorithm. When the quality of the result is the ultimate goal (regardless of the computation time), one can use the SCL algorithm after the SOM one to refine the solution (or in other words, one can increase the number of iterations performed without neighbors in the SOM procedure). Indeed even if a classical SOM usually ends without neighbor, the number of iterations performed without neighbors could be not sufficient to reach an optimal solution. Performing a sufficient number of iterations both with and (then) without neighbors is thus important for the quality of the solution.

## 6. Conclusion

The experiments illustrated in this paper, as well as many other ones performed on other data sets, indicate that the quality of the SOM algorithm resides not only in its topology preservation property, but also its vector quantization one. The SOM algorithm may be recommended compared to other VQ techniques like SCL, in order to reach a better minimum of the distortion error with a fixed number of iterations, or to reach faster a similar value of the distortion.

The better convergence properties cannot be proven theoretically. Nevertheless, we can make the analogy with simulated annealing techniques: the use of a neighborhood in the SOM algorithm introduces apparent disorder, making it possible to escape from a local minimum of the objective function and to increase the slope of convergence. Ending the VQ procedure with the SCL algorithm may be compared to ending a simulated annealing technique with a "temperature" parameter equal to zero.

The fact that the SOM method converges faster than classical VQ algorithms is well known; however, few studies bring systematic evidences supporting this claim. This paper shows quantitatively, through extended simulations, that the speed of convergence is increased by the use of the neighborhood property in the SOM. On artificial databases, where the optimal VQ solution may be computed analytically, we measured the VQ speed by the difference between the values of the quantizers during the convergence and their optimal stable solutions. On real databases where the optimal VQ solutions cannot be computed, we measured the VQ speed by the decrease of the distortion error with respect to the number of iterations. In both cases, we showed quantitatively how the use of the SOM algorithm, with fixed or decreasing number of neighbors, accelerates the convergence of the VQ learning. This produces clear indications that SOMs should be preferred to conventional VQ methods, even if topology preservation is not looked for.

## Acknowledgements

Michel Verleysen is a Senior Research Fellow of the Belgian National Fund for Scientific Research (FNRS). The authors thank the reviewers for their comments and suggestions.